\documentclass[10pt,a4paper]{amsart}
\usepackage{amssymb,amsfonts,amsthm,amsmath}
\usepackage{graphicx} 

\theoremstyle{plain}

\numberwithin{equation}{section} \numberwithin{theorem}{section}
\numberwithin{lemma}{section} \numberwithin{definition}{section}
\numberwithin{corollary}{section} \textheight =24cm
\begin{document}
\title[Magic Hexagons]{Magic Hexagon Formulas}
\author{Geoffrey B Campbell}
\address{Mathematical Sciences Institute \\
         The Australian National University \\
         ACT, 0200, Australia}

\email{Geoffrey.Campbell@anu.edu.au}

\begin{abstract}
We give a variety of magic hexagons of Orders from 3 to 7, many of which are extensions of known results.
We also give a theorem that their are an infinite number of magic hexagons of Order $n$ for any fixed positive integer $n$
for any arbitrary magic sum $M=m$ to be any desired integer $m$. We instigate theory and ideas associated with formula-based versions of the magic hexagons, which seem to be new.
\end{abstract}
\maketitle
\section{Magic Hexagons of Orders $3$, $4$, and $5$} \label{intro}

A magic hexagon of Order $n$ is an arrangement of numbers in a centred hexagonal pattern with $n$ ~cells on each edge, where the numbers in each row, and each diagonal, add to the same magic sum $M$. A normal magic hexagon contains the consecutive integers from $1$ to $3n^2-3n+1$. Normal magic hexagons exist only for $n=1$ with just one cell, and $n=3$ with $19$ cells, as shown below with $M=38$. There are also 26 distinct Order 3 magic hexagons with $M=0$. (see Wikipedia \cite{Wiki2025})

 \begin{figure} [!htbp]
\centering
    \includegraphics[width=9.61cm,angle=0,height=4.00cm]{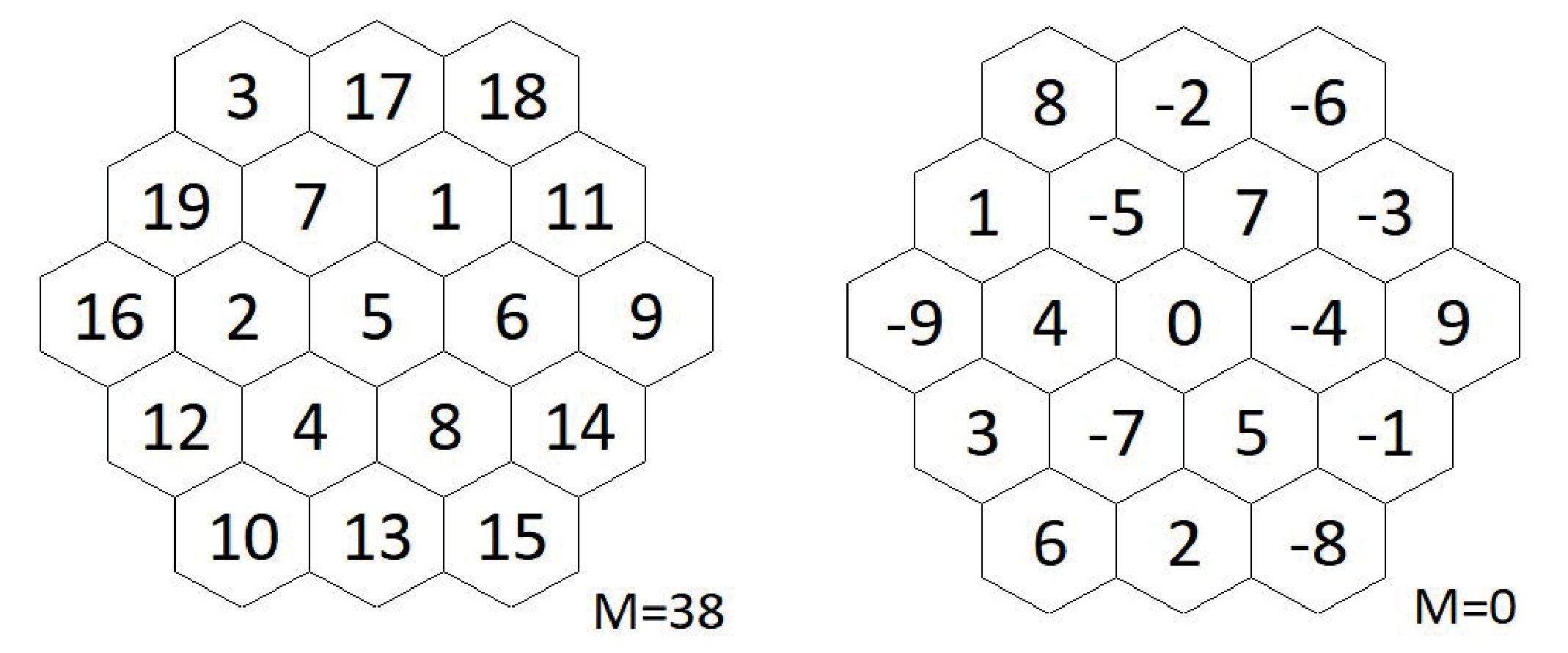}
  \caption[Order 3 Magic Hexagons: Normal (left) and Sum Zero (right).]{Order 3 Magic Hexagons: Normal (left) and Sum Zero (right).} \label{Fig1}
\end{figure}

The Order 3 normal hexagon was first published in 1887 and then rediscovered several times in the 20th century. For Orders $4,5,6,7,8$ and $9$, there are known abnormal and/or zero sum variants mostly from the 21st century. A known Order 4 abnormal hexagon with $M=111$ has entries 3 to 39, and a known Order 5 example with $M=244$ has entries 6 to 66 as shown here.

\begin{figure} [!htbp]
\centering
    \includegraphics[width=10.10cm,angle=0,height=4.00cm]{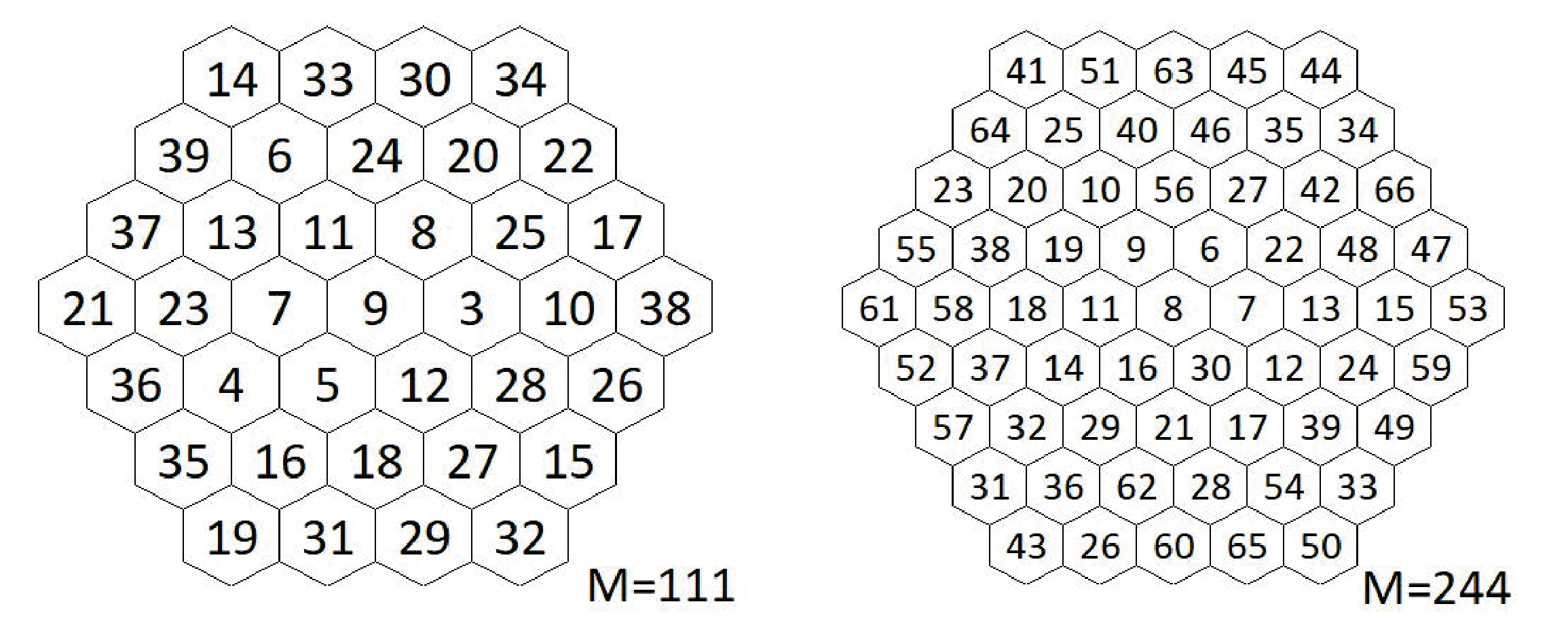}
  \caption[Order 4 and Order 5 Abnormal Magic Hexagons.]{Order 4 and Order 5 Abnormal Magic Hexagons.} \label{Fig2}
\end{figure}

\section{Formula-based Magic Hexagons of Orders $3$, and $4$.} \label{Orders34}

\begin{figure} [!htbp]
\centering
    \includegraphics[width=12.00cm,angle=0,height=5.00cm]{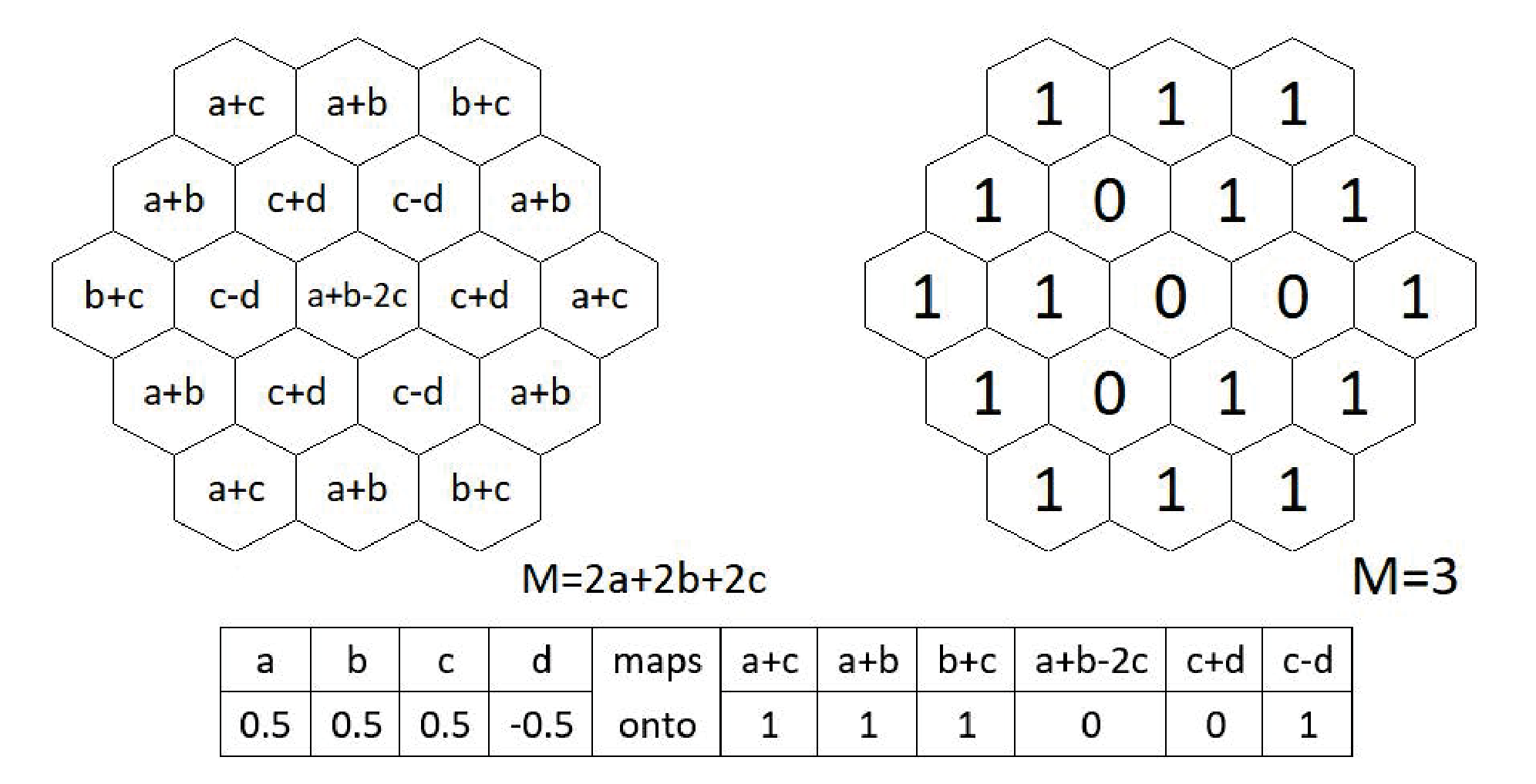}
   \caption[Order 3 Magic Hexagon formula with an $M=3$ case.] {Order 3 Magic Hexagon formula with an $M=3$ case.} \label{Fig3}
\end{figure}

We give formula-based, Order 3 and Order 4 magic hexagons, allowing duplicated number entries. Figure ~\ref{Fig3} has an Order ~3 version for $a$, $b$, $c$, and $d$ real numbers, with $M=2a+2b+2c$ (left), and a case with $M=3$ using only $1$s and $0$s (right). These formula-based magic hexagons first appeared in the author's small article recently at Campbell ~\cite{Campbell2025}.

The mappings to the $M=3$ magic hexagon are shown in the figure \ref{Fig3} table. Applying the left side hexagon of figure \ref{Fig1}, we can add to corresponding cell multiples of the right hexagon in figure ~\ref{Fig5}, to obtain two Order 3 magic hexagons, each with two sequences of consecutive integers. Figure ~\ref{Fig6} contains mostly squares and also derives from LHS hexagon in figure ~\ref{Fig3}.

In figure ~\ref{Fig3}:

 \begin{enumerate}
   \item the LHS hexagon has $M=59$ and values in the two sequences: $4,5,\cdots,10$ and $15,16,\cdots,26$.
   \item the RHS hexagon has $M=10$ and values in the two sequences: $-7,-6,\cdots,8$ and $13,14,15$.
 \end{enumerate}

\begin{figure} [!htbp]
\centering
    \includegraphics[width=13.13cm,angle=0,height=4.00cm]{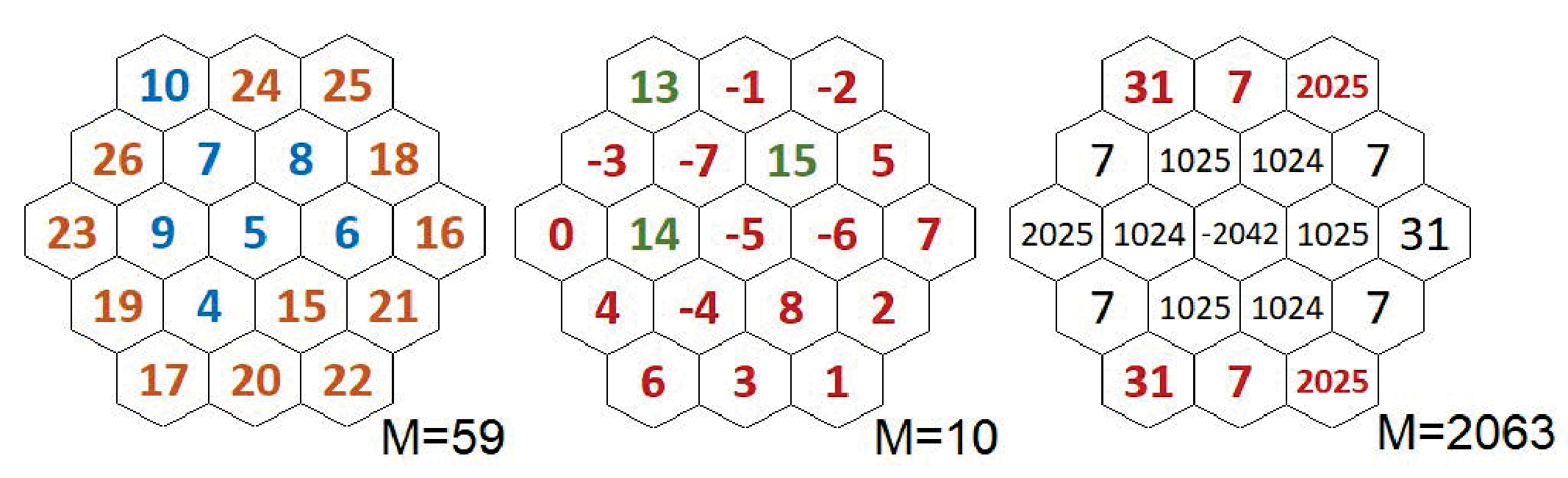}
  \caption[Order 3: Hexagons with consecutive numbers and 31/7/2025.] {Order 3: Hexagons with consecutive numbers and 31/7/2025.} \label{Fig4}
\end{figure}

\begin{figure} [!htbp]
\centering
    \includegraphics[width=12.94cm,angle=0,height=4.00cm]{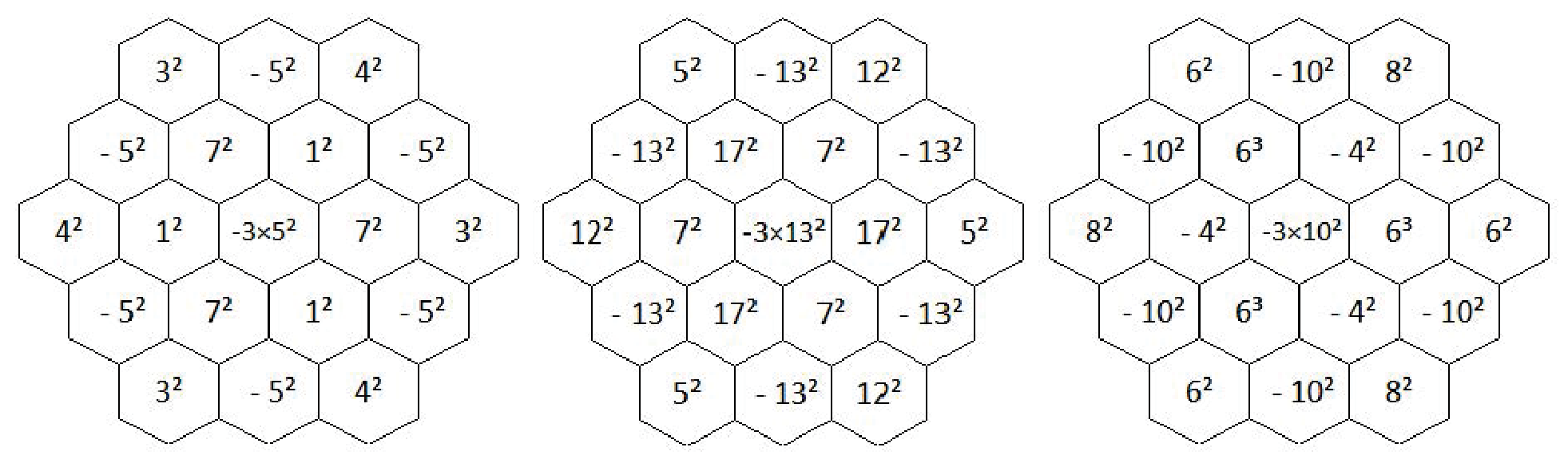}
   \caption[Order 3 Magic Hexagons each M=0 with mostly square entries.] {Order 3 Magic Hexagons each M=0 with mostly square entries.} \label{Fig5}
\end{figure}

 Alternatively, we can contrive values by solving the simultaneous equations as in the AU format date 31st July 2025 of figure ~\ref{Fig4}. ie. We solved here $a+c=31$, $a+b=7$, $b+c=2025$, $d=\frac{1}{2}$.

The three hexagons in figure ~\ref{Fig5} were found from considering cases $a=A$, $b=B^2$, $c=-\;C^2$ and then substituting Pythagorean triples.

\newpage
Also, we can derive an Order 4 formula-base example (see figure ~\ref{Fig6}) with a case $M=83$ with 31/7/2025 AU Format Date:

 \begin{figure} [!htbp]
\centering
    \includegraphics[width=10.27cm,angle=0,height=5.00cm]{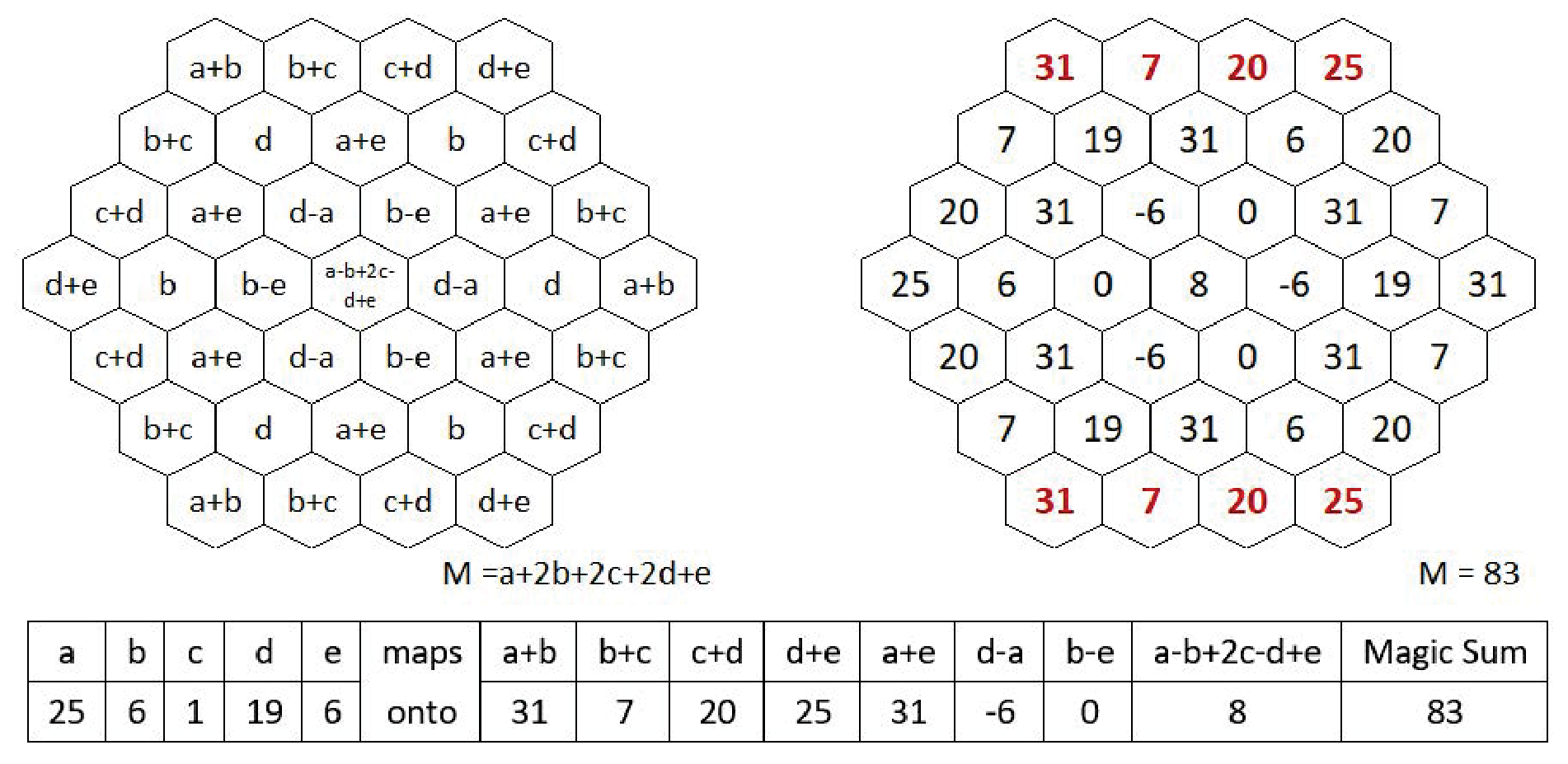}
  \caption[Order 4: Formula-based Magic Hexagon with a case.]{Order 4: Formula-based Magic Hexagon with a case.} \label{Fig6}
\end{figure}

\begin{figure} [!htbp]
\centering
    \includegraphics[width=8.14cm,angle=0,height=5.00cm]{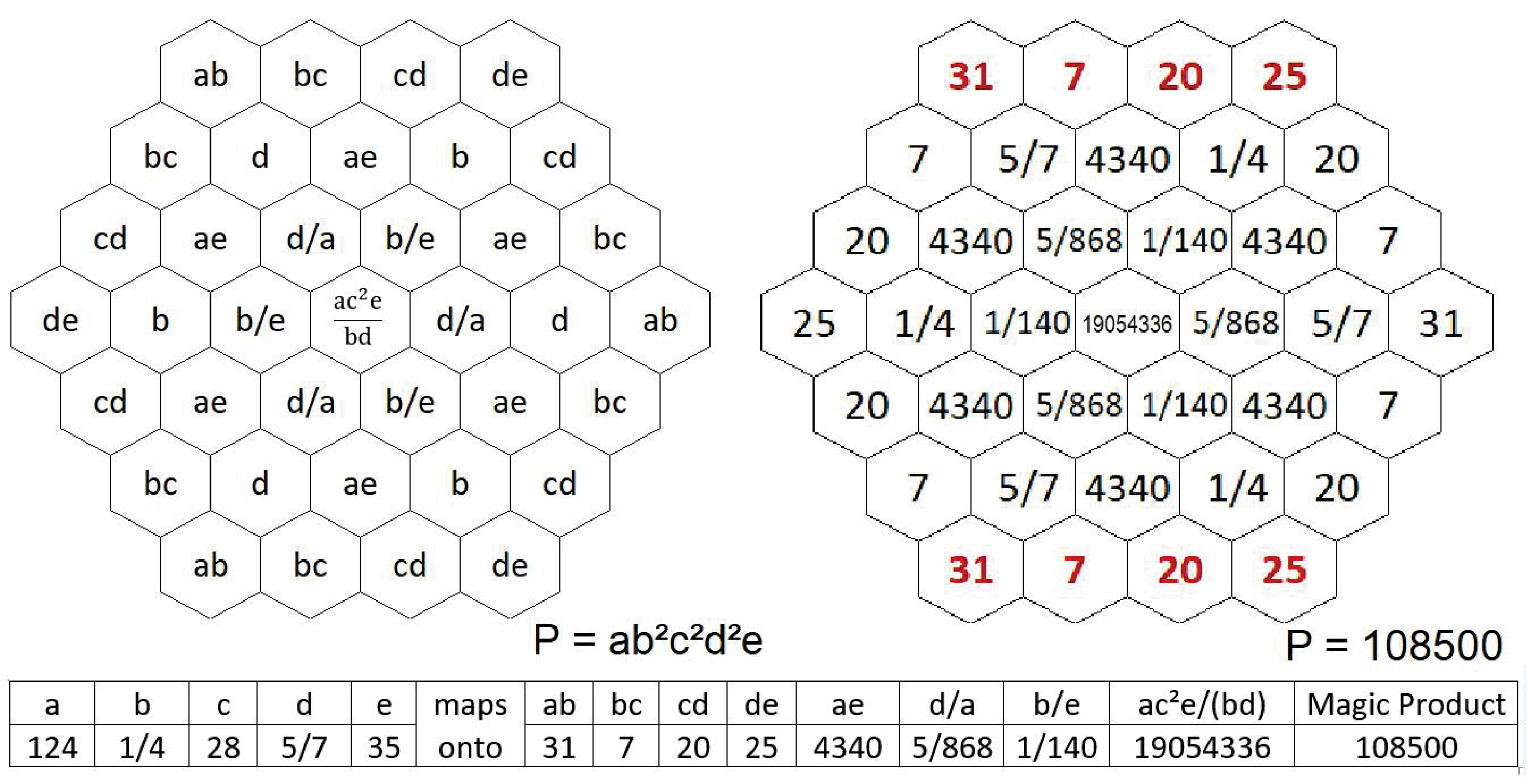}
  \caption[Order 4 Magic Hexagon Formula Products with 31/7/2025 case.] {Order 4 Magic Hexagon Formula Products with 31/7/2025 case.} \label{Fig7}
\end{figure}

Our formula-based magic hexagons seem to be new. It is easy to find them for Orders $5$ and $6$, which we do in the next section. The original 1887 normal hexagon proof (von Haselberg \cite{vonHaselberg1887}) is in Bauch's 2012 paper  \cite{Bauch2012}. Our modern figure ~\ref{Fig2} abnormal magic hexagons are in the Wikipedia article \cite{Wiki2025}, credited there to \textit{Arsen Zahray} without proper citation. Also, the additive formulas imply multiplicative magic hexagons, so the Magic Sum $M=a+2b+2c+2d+e$ in figure \ref{Fig6}, becomes Magic Product $P=ab^2c^2d^2e$ in figure ~\ref{Fig7}.

\newpage

\begin{figure} [!htbp]
\centering
    \includegraphics[width=8.67cm,angle=0,height=5.00cm]{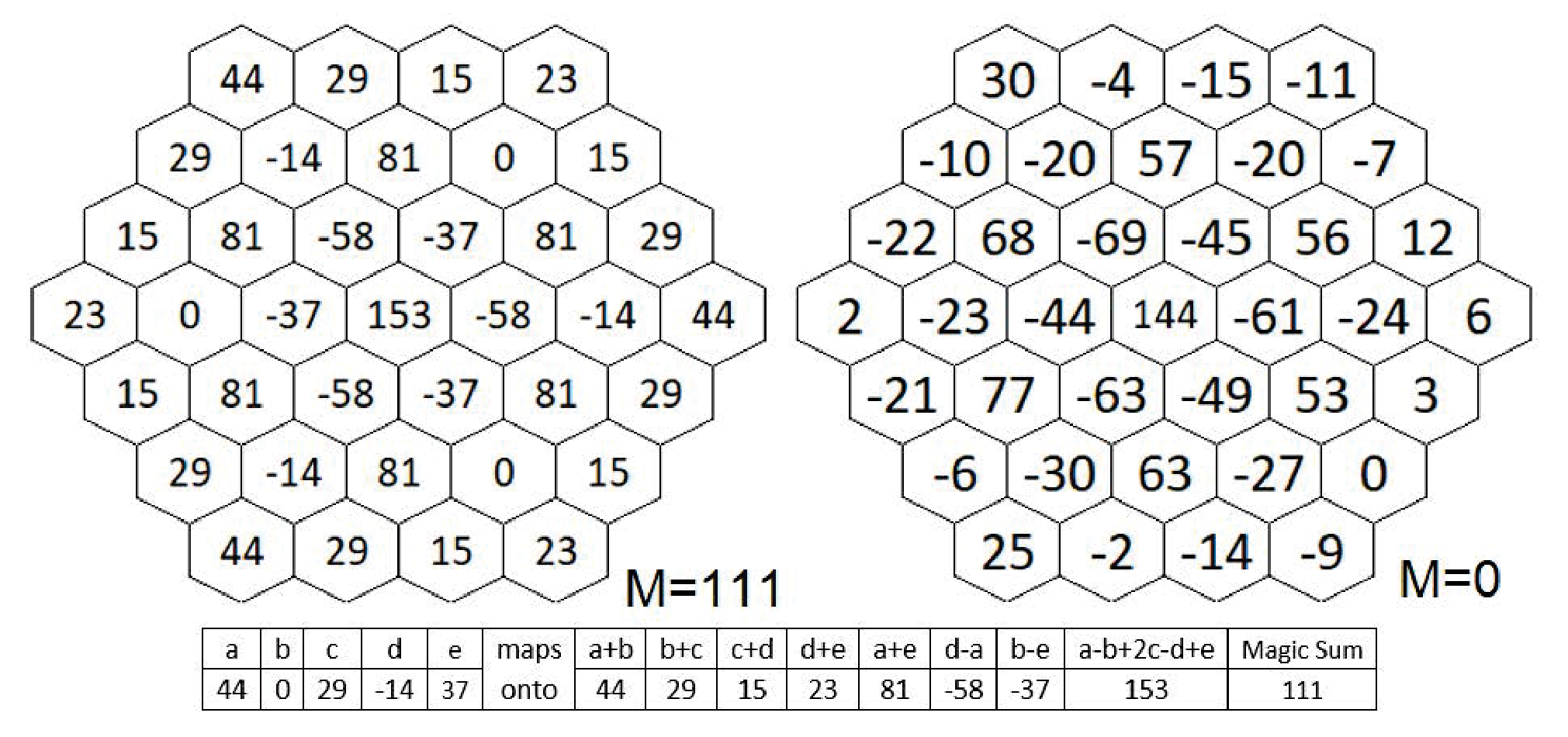}
  \caption[Formula-based Order 4 Magic Hexagon M=111 leading to an M=0 version.] {Formula-based Order 4 Magic Hexagon M=111 leading to an M=0 version.} \label{Fig8}
\end{figure}

From our figure \ref{Fig2}, Order 4 hexagon M=111 we can add multiples of formula-based Order 4 hexagons to obtain new magic hexagons. So, for example, taking $a$=44, $b$=0, $c$=29, $d$=-14, $e$=37 in figure \ref{Fig2}, we have the LHS hexagon in figure \ref{Fig8}. Subtracting this from LHS figure \ref{Fig2} hexagon cell-by-cell gives the RHS figure ~\ref{Fig8} hexagon.

\textbf{Sequential workings for Order 4 RHS figure ~\ref{Fig8} Hexagon M=0:}

The distinct numbers (apart from repetition of -20):
\begin{center}
$-69, -63, -61, -49, -45, -44, -30, -27, -24, -23, -22, -21, -20, -20, -15, -14,$

$ -11, -10, -9, -7, -6, -4, -2, 0, 2, 3, 6, 12, 25, 30, 53, 56, 57, 63, 68, 77, 144,$
\end{center}
appear in the Order ~4 Magic Hexagon with M=0,

\begin{figure} [!htbp]
\centering
    \includegraphics[width=10.46cm,angle=0,height=4.00cm]{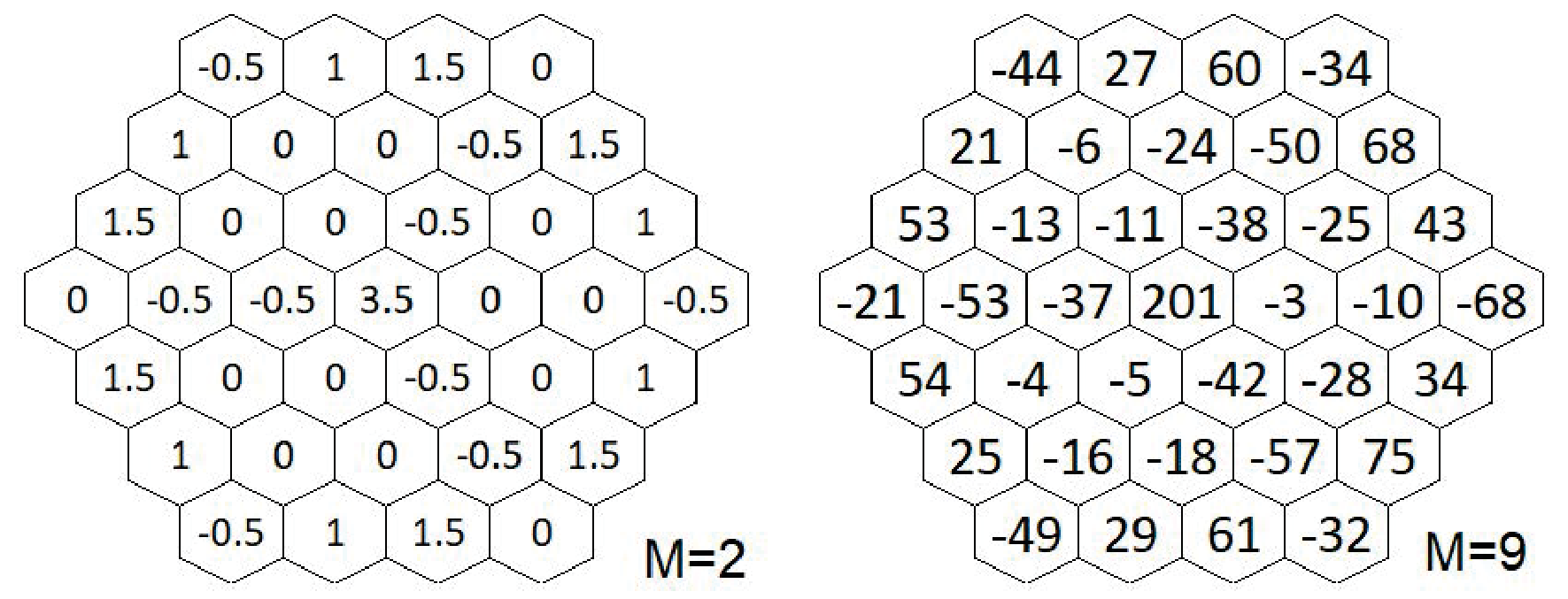}
  \caption[Formula-based Order 4 Magic Hexagon M=2 leading to an M=9 version.] {Formula-based Order 4 Magic Hexagon M=2 leading to an M=9 version.} \label{Fig9}
\end{figure}

\textbf{Sequential workings for Order 4 Hexagon M=9:} Next we give an example of a simple M=2 Order 4 magic hexagon which forms an M=9 magic hexagon from a linear combination of $-A_1+60A_2=A_3$; where $A_1$ is the LHS hexagon in figure ~\ref{Fig2}, and $A_2$ is the LHS hexagon in figure ~\ref{Fig9}, with $A_3$ the RHS hexagon in figure ~\ref{Fig9}. The RHS hexagon of figure ~\ref{Fig9} is the Order 4, $M$=9 hexagon $A_3$ having the distinct entries listed here in numerical order:

\begin{center}
$-68, -57, -53, -50, -49, -44, -42, -38, -37, -34, -32, -28, -25, -24, -21, -18, -16,$

$-13, -11, -10, -6, -5, -4, -3, 21, 25, 27, 29, 34, 43, 53, 54, 60, 61, 68, 75, 201.$
\end{center}

So the RHS figure ~\ref{Fig9} hexagon is derived from the LHS M=2 hexagon with fractional entries.

\begin{figure} [!htbp]
\centering
    \includegraphics[width=15.00cm,angle=0,height=7.84cm]{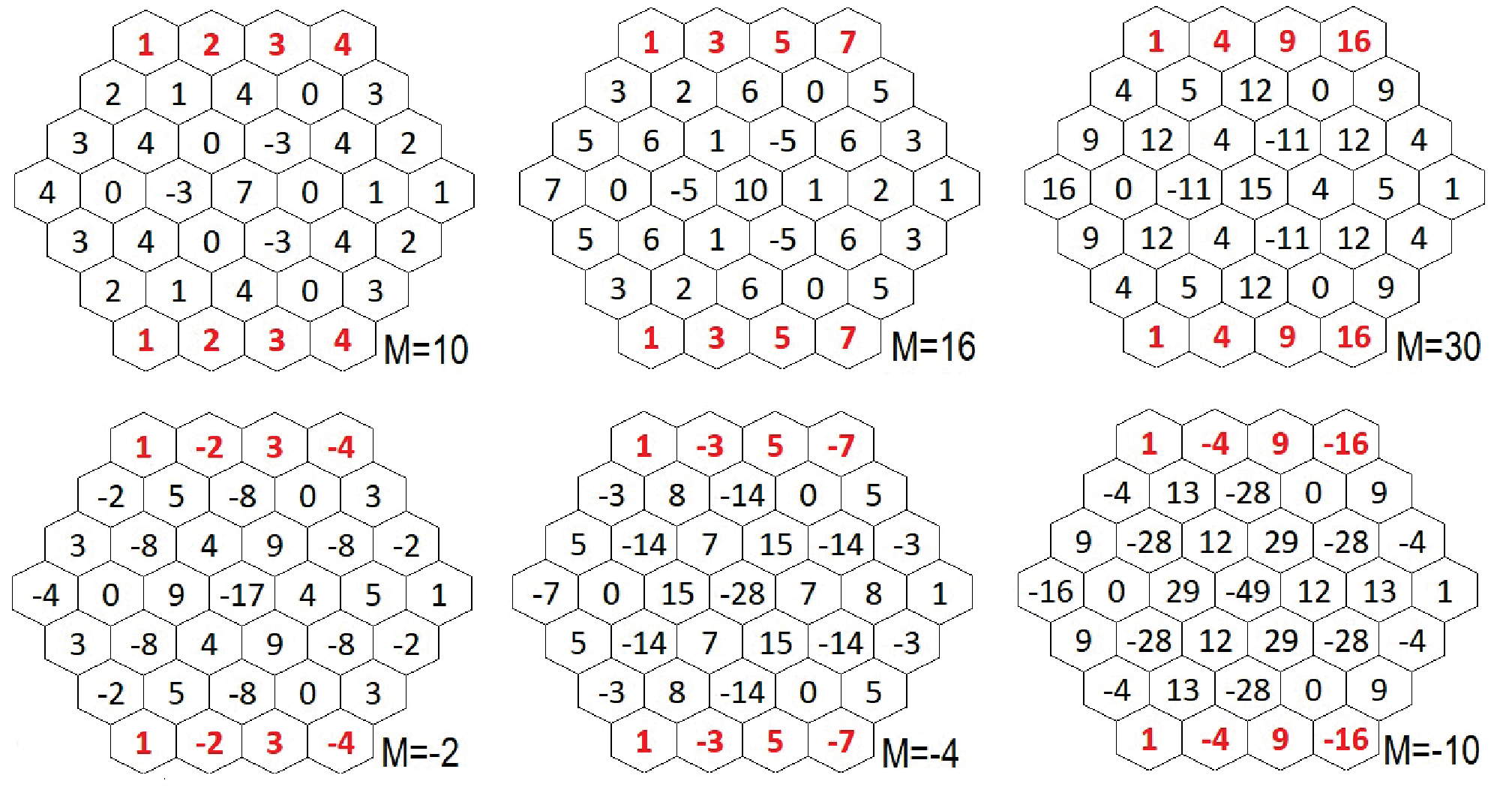}
  \caption[Six Derived Order 4 Magic Hexagons: M=10, 16, 30, -2, -4, -10.] {Six Derived Order 4 Magic Hexagons: M=10, 16, 30, -2, -4, -10.} \label{Fig10}
\end{figure}

In figure ~\ref{Fig10} we give further examples of derived Order 4 magic hexagons

Clearly the variety of possible combinations to obtain a particular value of $M$ is infinite.
\newpage

\section{Formula-based Magic Hexagons of Orders $5$, and $6$.} \label{Orders56}

\begin{figure} [!htbp]
\centering
    \includegraphics[width=12.32cm,angle=0,height=6.50cm]{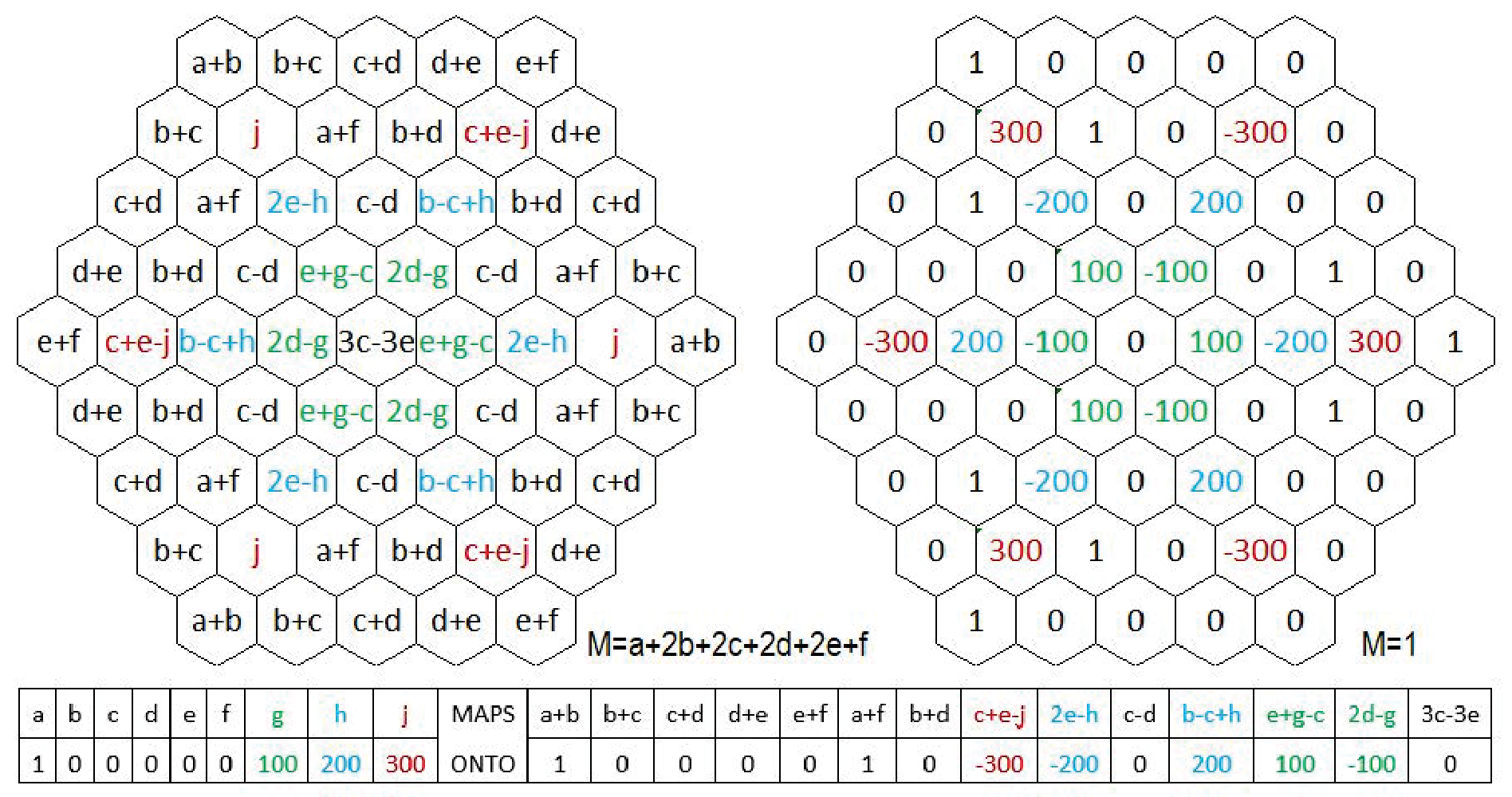}
  \caption[Order 5 Formula-based M=a+2b+2c+2d+2e+f mapped onto an M=1 hexagon.] {Order 5 Formula-based M=a+2b+2c+2d+2e+f mapped onto an M=1 hexagon.} \label{Fig11}
\end{figure}

We state a formula-based Order 5 magic hexagon. This is easy to derive by working from the outside layer of hexagons inward, and maintaining symmetry wherever possible. The figure \ref{Fig11} formula entries are colour-coded and the example for $M=1$ of Order 5 magic hexagons demonstrates how to construct an infinite set of $M=1$ solutions.

\begin{figure} [!htbp]
\centering
    \includegraphics[width=10.30cm,angle=0,height=8.00cm]{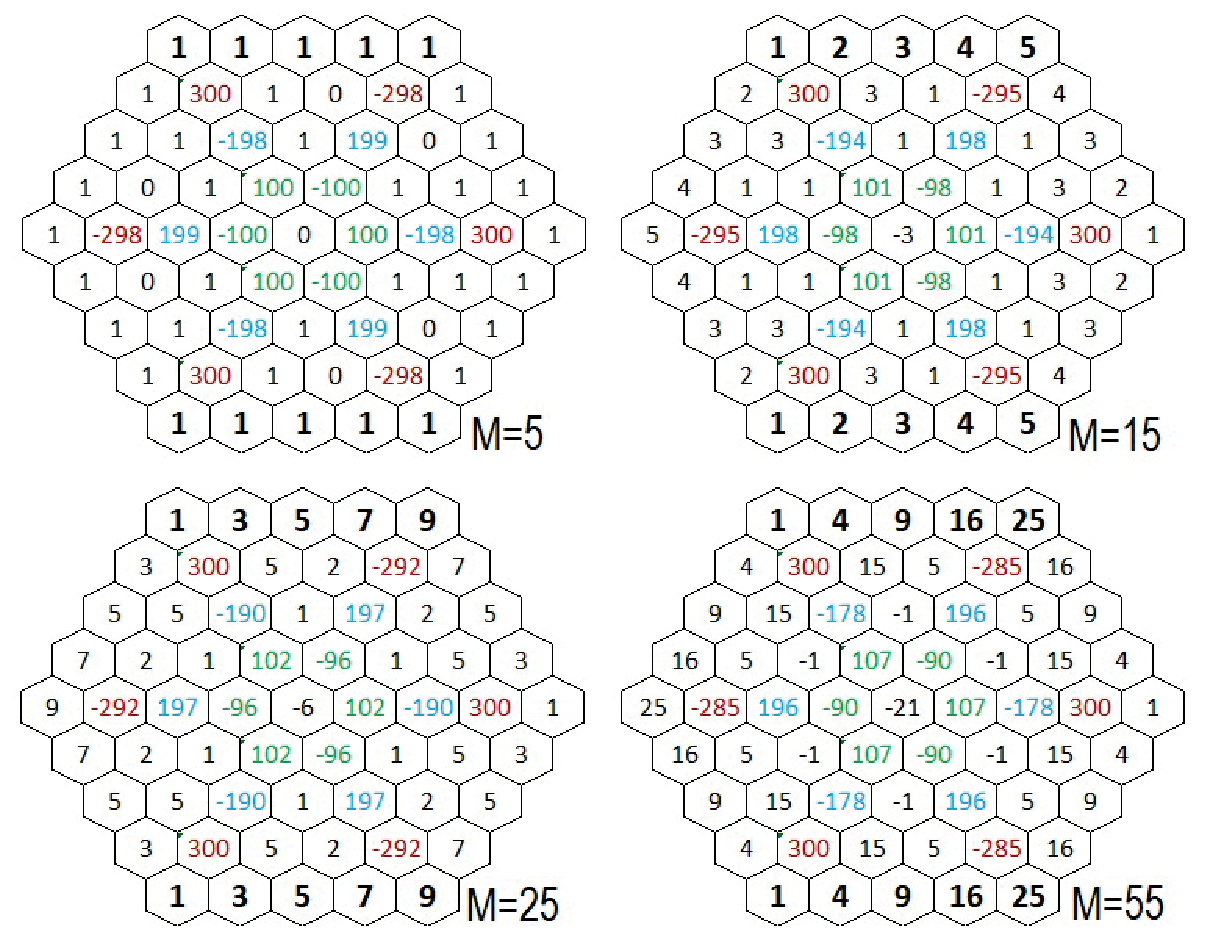}
  \caption[Order 5 derived Magic Hexagon cases: M=5, 15, 25, 55.] {Order 5 derived Magic Hexagon cases: M=5, 15, 25, 55.} \label{Fig12}
\end{figure}

We give an Order 6 formula-based magic hexagon, mapped to a case with every outer cell value 91 contriving an M=546 magic hexagon in figure \ref{Fig11}. We then subtract the cell-by-cell values of this hexagon from the known M=546 magic hexagon created by Louis Hoelbling, October 11, 2004 (see Wikipedia \cite{Wiki2025}) which we reproduce here in LHS of figure \ref{Fig12}, to obtain an M=0 version on RHS of figure \ref{Fig12}.

 \begin{figure} [!htbp]
\centering
    \includegraphics[width=14.67cm,angle=0,height=8.00cm]{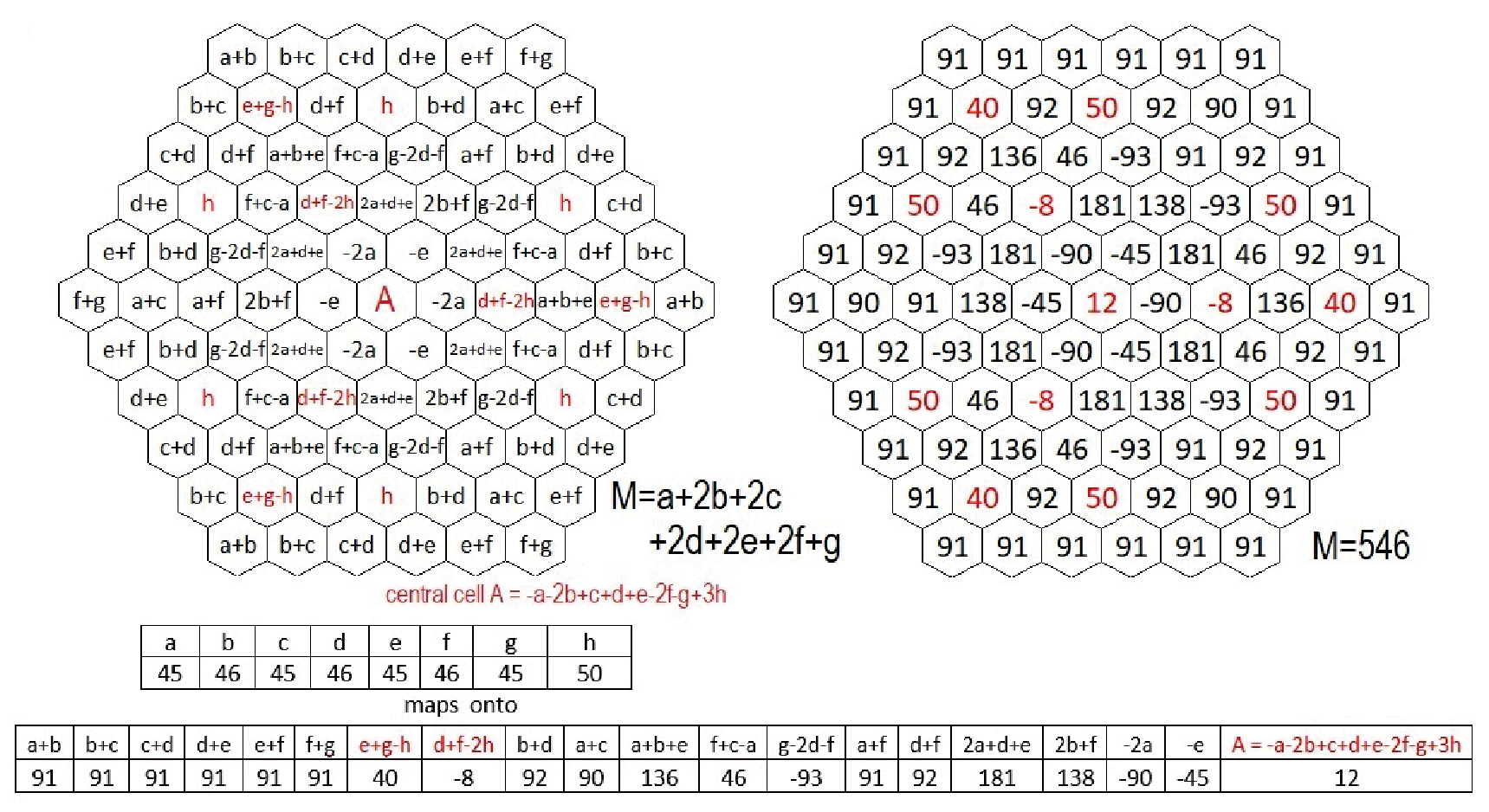}
  \caption[Order 6 Formula-based M=a+2b+2c+2d+2e+2f+g mapped onto an M=546 hexagon.] {Order 6 Formula-based M=a+2b+2c+2d+2e+2f+g mapped onto an M=546 hexagon.} \label{Fig13}
\end{figure}

 \begin{figure} [!htbp]
\centering
    \includegraphics[width=14.44cm,angle=0,height=6.00cm]{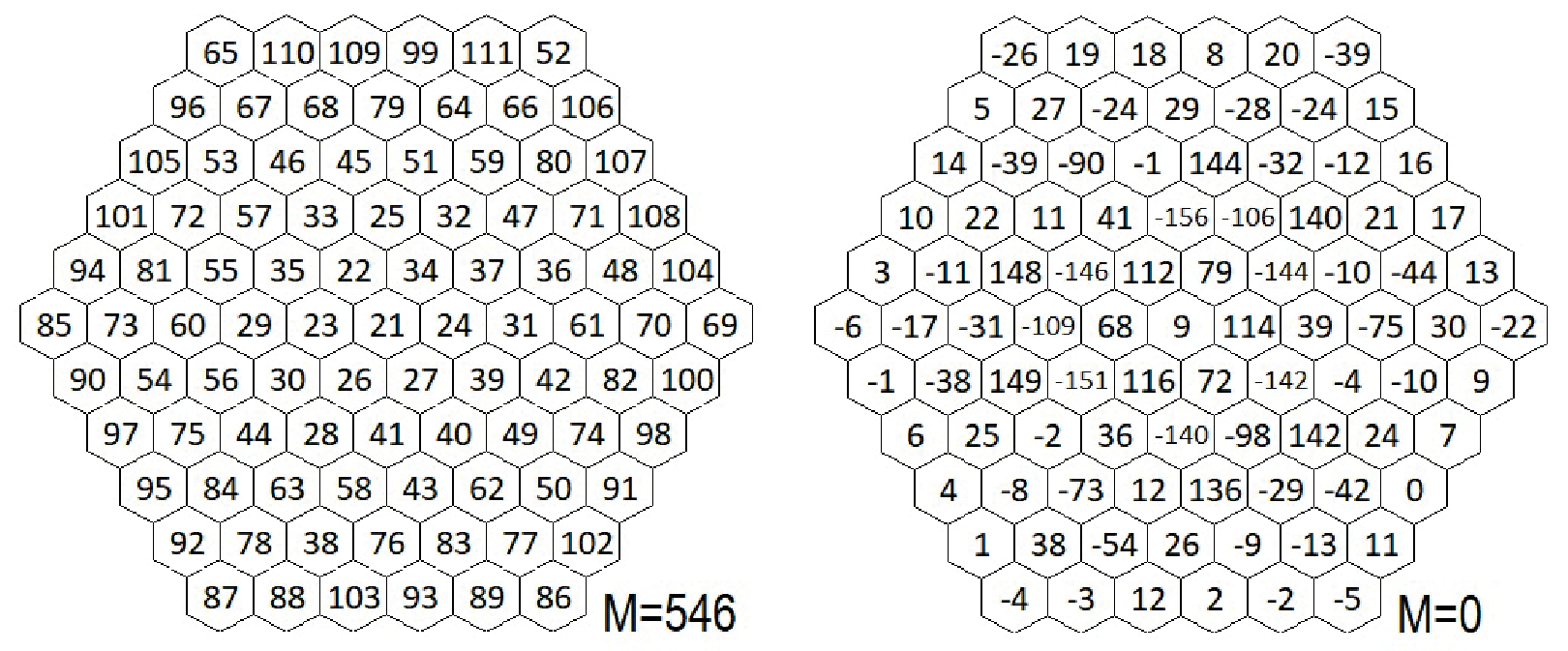}
  \caption[Order 6: M=546 hexagon 21 to 119 entries; and derived M=0 hexagon.] {Order 6: known M=546 hexagon 21 to 119 entries; and derived M=0 hexagon.} \label{Fig14}
\end{figure}

\section{Order 7 Magic Hexagons and beyond.} \label{Orders7}

\begin{figure} [!htbp]
\centering
    \includegraphics[width=15.00cm,angle=0,height=14.33cm]{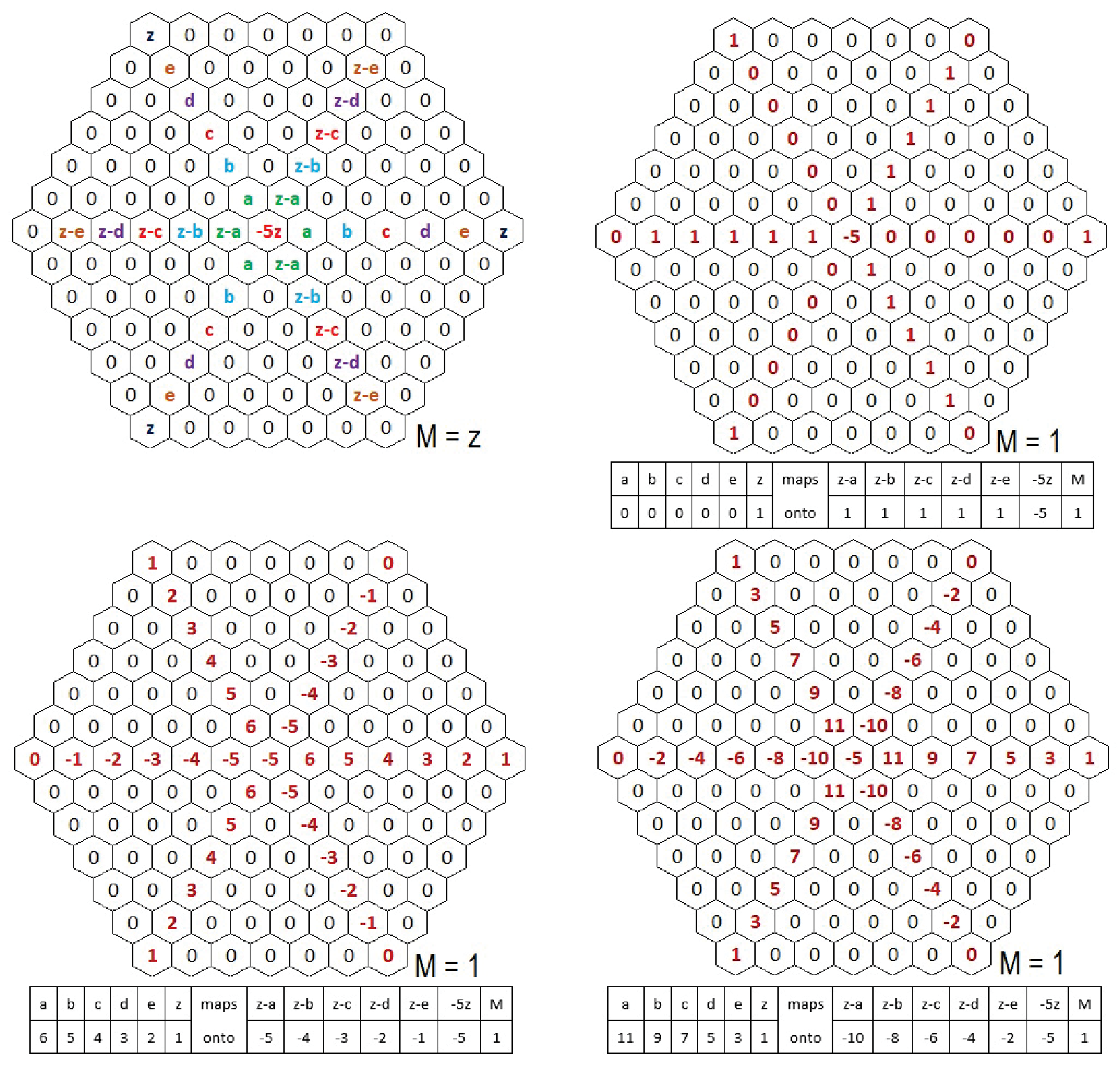}
  \caption[Order 7 Formula-based M=z mapped onto three simple M=1 hexagons.] {Order 7 Formula-based M=z mapped onto three simple M=1 hexagons.} \label{Fig15}
\end{figure}

We state a formula-based Order 7 magic hexagon with M=z with z an arbitrary integer, and derive three cases with $M=1$ solutions in figure \ref{Fig15}. This then suggests, and this is easy to prove in generalization, that the diagonal method of figure \ref{Fig15} is applicable to magic hexagons of any order and for any desired value of the magic sum M.

 \begin{figure} [!htbp]
\centering
    \includegraphics[width=12.56cm,angle=0,height=12.00cm]{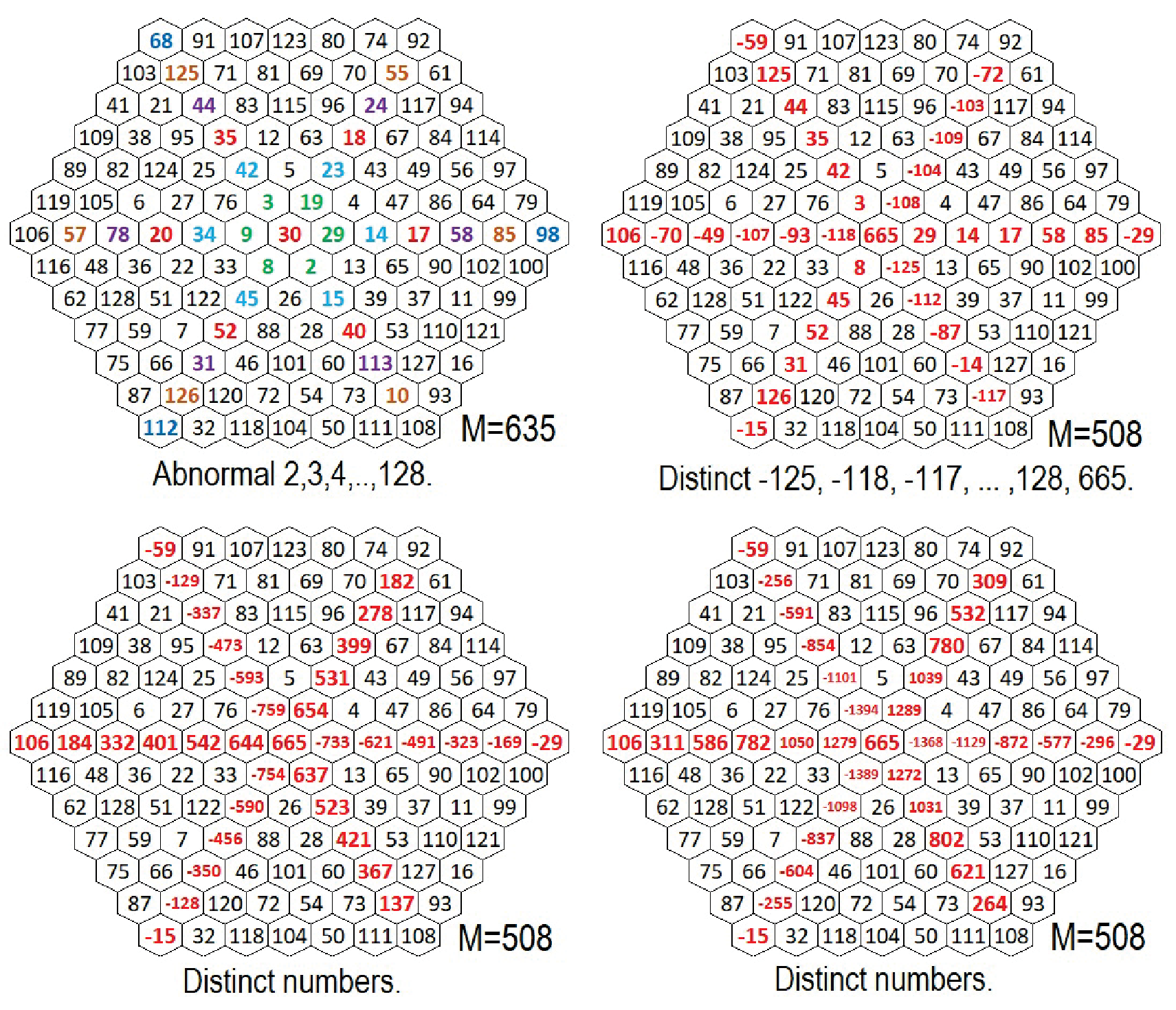}
  \caption[Order 7 known M=635 abnormal mapped onto three M=508 hexagons.] {Order 7 known M=635 abnormal mapped onto three M=508 hexagons.} \label{Fig16}
\end{figure}

 \begin{figure} [!htbp]
\centering
    \includegraphics[width=15.00cm,angle=0,height=9.83cm]{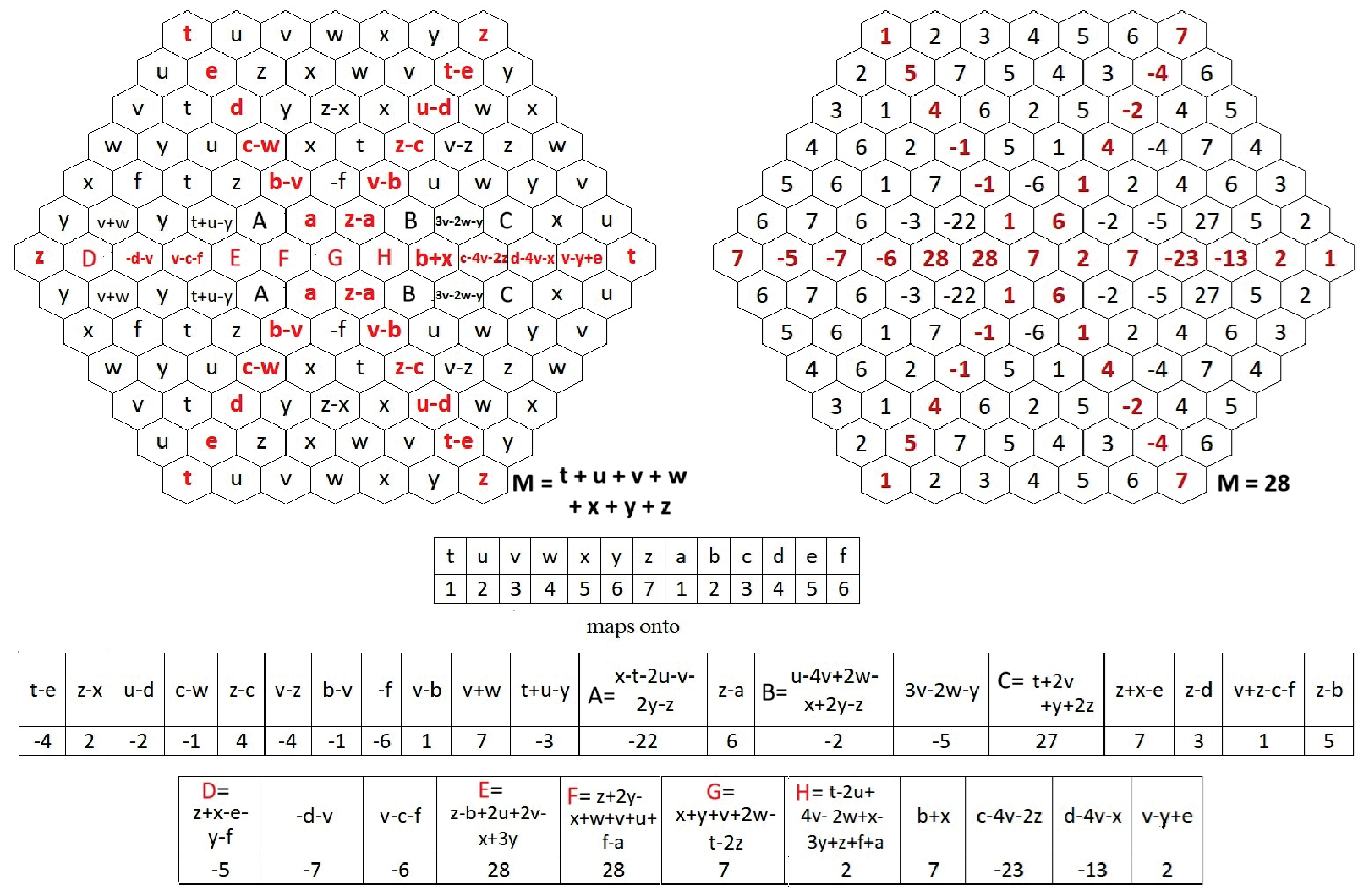}
  \caption[Order 7 M=t+u+v+w+x+y+z formula with M=28 case.] {Order 7 M=t+u+v+w+x+y+z formula with M=28 case.}
  \label{Fig17}
\end{figure}


\newpage







\begin{thebibliography}{99}
\bibitem{Bauch2012}
Bauch, H. ``125 Jahre magisches Sechseck." Internat. Math. Nachrichten. Nr. 219 (2012), 13–24. (https://www.oemg.ac.at/IMN/imn219.pdf).
\bibitem{Campbell2025}
Campbell, G. B. ``Fun with numbers 16: Some Magic Hexagons'', Aust. Math. Soc. Gazette, Volume 52, No3, July 2025. (https://austms.org.au/wp-content/uploads/2025/08/523Web\_4.pdf)
\bibitem{Trigg1964}
Trigg, C. W. ``A Unique Magic Hexagon." Recr. Math. Mag. 46, 40-43, Jan./Feb. 1964. (http://www.mathematik.uni-bielefeld.de/~sillke/PUZZLES/magic-hexagon-trigg).
\bibitem{Trigg1972}
Trigg, C. W. ``P824: A Well-Known Magic Hexagon." Math. Mag. 45, 100, 1972.
\bibitem{Trigg1973}
Trigg, C. W. ``Solution to Problem P824." Math. Mag. 46, 44-45, 1973.
\bibitem{vonHaselberg1887}
von Haselberg, E. ``Problem and Solution of the Unique Magic Hexagon of Order 3." Manuscript, 1887.
\bibitem{Weisstein2025}
Weisstein, E. W. ``Magic Hexagon." From MathWorld--A Wolfram Web Resource, viewed June 2025. (https://mathworld.wolfram.com/MagicHexagon.html)
\bibitem{Wiki2025}
Wikipedia, ``Magic hexagon.", the free encyclopedia, viewed June 2025. (https://en.wikipedia.org/wiki/Magic\_hexagon)
\end{thebibliography}
\end{document}